\documentclass[11pt]{article}
\usepackage{amstext,amsmath,amssymb,amsthm}
\usepackage{latexsym}
\usepackage{exscale}

\newcommand{\ov}{\overline}

\newcommand{\dsize}{\displaystyle}

\newcommand{\R}{{ \mathbb  R  }}

\newcommand{\N}{ \mathbb  N }
\newcommand{\Z}{ \mathbb  Z }
\renewcommand{\L}{{\mathcal  L} }
\newcommand{\B}{{\mathcal  B} }

\newcommand{\BEL}{\begin{equation}\label}
\newcommand{\EE}{\end{equation}}
\renewcommand{\t}{\tilde}

\renewcommand{\medskip}{\vskip .5 cm}
\newtheorem{Thm}{Theorem}[section]
\newtheorem{Lemma}[Thm]{Lemma}

\newtheorem{Prop}[Thm]{Proposition}

\begin{document}

 \centerline{\bf \Large  Exponential bases    on two dimensional  trapezoids }

\bigskip
 
 \centerline{\bf   Laura De Carli  and Anudeep Kumar}

  \medskip
 \begin{abstract}
We discuss   existence and    stability of Riesz bases   of exponential type of $L^2(T)$  for special domains $T\subset\R^2$ called trapezoids. 
  We construct  exponential   bases  on   $L^2(T)$ when $T$ is a    finite  union  of   rectangles  with the same height. We also generalize our main  theorems in dimension $d\ge 3$.

 \medskip\noindent Mathematics Subject Classification:   42C15,  42C30.

\noindent
{Keywords:  exponential bases, trapezoids, multirectangles.}  \end{abstract}

 \section{Introduction}

 \setcounter{section}{1} \setcounter{Thm}{0} \setcounter{equation}{0}

 The study of Riesz bases   in Hilbert spaces  is a fruitful topic of investigation since many decades.    Riesz  bases allow to decompose   $L^2$ functions  (which can be thought of as  signals) in a unique way into a sum of  basic frequencies.  
Riesz bases  of exponential type,   i.e.,    in the form of 
${\cal B}=\{e^{i\alpha_n(x)}\}_{n\in\Z}$, with $\alpha_n\in L^\infty (D)$,  are   especially simple and easy to work with.   Here, $D$ is a domain of $\R^d$, i.e., a bounded and measurable set of finite  Lebesgue measure.  
If the
$\alpha_n(x)$'s are linear, ${\cal B}$  is called   {\it exponential basis}.
Exponential bases  of $L^2(a,b)$ are well   studied since   when 
Paley and Wiener \cite{PW} and Levinson  \cite{L}   explored the  possibility of  non-harmonic Fourier series (see also \cite{B}, \cite{DS}). Since then,   exponential  bases  have been used in various fields of mathematics such as the theory of non-selfadjoint operators, control theory   and signal processing. See   \cite{AI}, \cite{Be}, \cite{FP}, just to cite a few.  
%
  Much less is  known  about the existence of exponential bases of $L^2(D)$   when  $D$ has dimension $d\ge 2$.  In \cite{LR}  the existence of 
 exponential bases    is proved  for  convex 
symmetric polygon.  To the best of our knowledge,  the only results in the literature  in dimension greater that two are in   \cite{M} and in the very recent   \cite{GL}.  

Exponential bases  in dimension $d=2$ and $d=  3$   are    important for the applications to  image recognition,  magnetic resonance imaging and data compression.  We refer the interested reader to  the recent \cite{CK} and   the references cited there. 

In this paper we   construct Riesz bases of exponential type  on a class of   domains  of $\R^2$ that, with some abuse of terminology, are called {\it trapezoids}. 
Let  $-\infty<a<b<\infty$, and let $f:[a,b]\to\R$.  Assume that      \begin{equation}\label{e-admis} 0< \ell <f(x)< M \quad \mbox{   almost everywhere (a.e) on $[a,b]$}.
\end{equation} We refer to the set  
\begin{equation}\label{e-Trapez} T=\{ (x,y)\in\R^2 : \ |x| \leq f(y),  \   a\leq y \leq b \}  \end{equation} 
 as to the   {\it trapezoid    bounded by $f$}.  
 
In what follows  we will   assume  that  trapezoids are bounded   by   functions that satisfy \eqref{e-admis}, and we will   let  $[a,b]= [0,1]$ for the sake of simplicity.

A trapezoid  bounded  by a step function  is    union of a finite family of rectangles. These trapezoids  are especially interesting
because   step functions are dense in  $L^p(0,1)$ for every $p<\infty$.  We say  that a positive  step function   is {\it regular} if it is    constant on  intervals of equal length;  thus, a regular step function on $[0, 1]$ can be written as \begin{equation}\label{eregular-step} s(x)= 
 \sum_{j=1}^{N-1} b_j \chi_{[\frac {j-1} N,\ \frac{j}N)}(x)  + b_N  \chi_{[\frac {N-1} N,\ 1]}(x)   
\end{equation}
where $\chi$ is the characteristic function, $ N$    a positive integer  and   $b_j >0$. A trapezoid  bounded by a    step function  as in \eqref{eregular-step} is called  {\it multi-rectangle  with $N$ steps}.

  In  Section 3 we approximate  the characteristic function of a trapezoid with 
 regular step functions to construct special Riesz bases of exponential type on $L^2(T)$. 

\begin{Thm}\label{T-Stepf}  
Let $T $ be the trapezoid bounded by  a piecewise continuous  $f$. 
There exists a family of regular step functions 
$\{s_n \}_{n\ge 1} $  which converge  uniformly   to $f$, for which
\begin{equation}\label{e-stepbase} {\mathcal B} =\left \{   e^{\pi i(\frac{ n x }{s_n(y)} +  2  ky) }\right\}_{n,k\in\Z} \end{equation} 
is a Riesz basis of  $L^2(T)$.
 \end{Thm}
Recall that $f$   is piecewise continuous  if it is continuous  everywhere except at
a finite number of points.
Theorem \ref{T-Stepf} allows to construct a  Riesz basis  of $L^2(T)$  in the form of  $\{e^{  i (x \lambda_n(y)+2k y)}\}_{n, k\in\Z}$, with   $\lambda_n(y) = n/s_n(y)$ piecewise constant.  If $\lambda_n$ was  constant, then ${\cal B}$ would be an exponential basis of $L^2(T)$,  
but the piecewise constant   functions   are the "next best thing" after the constants.

Theorem \ref{T-Stepf} is a corollary of  the following

\begin{Thm}\label{T-Main11}    Let $T$ be the trapezoid bounded by $f$. Let
  $\{g_n(y)  \}_{n\in\Z}  \subset L^\infty(0,1)$  be such that   
\begin{equation}\label{E-cond11} 
\left| \frac{ f(y)}{g_n(y)}-1\right|=\epsilon_n(y) < \frac 1{4|n| }  \quad \rm{a.e.}
 \end{equation}
  Then, 
 ${\mathcal B} =\left \{   e^{\pi i (\frac{ n x }{g_n(y)} +   2  ky) }\right\}_{n,k\in\Z} $
is a Riesz basis for $L^2(T)$. The constant $\frac 14$  in \eqref{E-cond11} cannot be replaced by any larger constant. 
 \end{Thm}

\noindent 
 This theorem can be viewed as an analogue of the celebrated  Kadec  stability  theorem, that we state and discuss in     Section  2.3.
In  Section 4 we show that  trapezoids bounded by  regular step functions have   exponential bases.  We prove the following

\begin{Thm}\label{T-Lev2}
Let    $R\subset [-1, 1]\times [0,1]$ be  a multi-rectangle with $N$ steps. There exists sequences  $\{n_k\}_{k\in\Z}$,  $\{m_h\}_{h\in\Z}\subset \Z$   such that  $\{e^{2\pi i (\frac{  n_k x}{2N}+ N m_h y)  }\}_{h, k\in\Z}$ 
is a  Riesz basis of $L^2(R)$.
 \end{Thm}

 To prove Theorem \ref{T-Lev2} we  show   that  an exponential basis of $L^2(R)$ can be obtained   from  an exponential basis on a finite unions of    disjoint segments of $\R$, (or   {\it multi-intervals}). 
 It has been recently proved by   G. Kozma  and  S. Nitzan   that every multi-interval has an exponential basis.  The proof of  Theorem \ref{T-Lev2} relies on  the main theorem in \cite{KN}.
 
  Interest in exponential Riesz bases on  multi-intervals has its roots
in practical applications to sampling of band-limited signals, and the first partial
results came from there.  The  main theorem  in \cite{KN} improves results of N. Lev \cite{L} and K. Seip \cite{S2}.

\noindent

%
The plan of this paper is the following.  In section 2 we   state some preliminary definitions and results. In Section 3 we prove  
Theorems \ref{T-Stepf} and \ref{T-Main11}. In Section 4 we prove Theorem \ref{T-Lev2}.  In Section 5 we prove analogues of our main theorems in dimensions $d\ge 2$.
We have stated some remarks and conjectures in Section 6. 

\medskip
\noindent
{\it Acknowledgement.} We wish to thank   K. Seip for   advice  and suggestions, and the referee for suggesting multidimensional versions of our main theorems.

\section{Preliminaries}
  \setcounter{Thm}{0} \setcounter{equation}{0}

In this section we  have collected some definitions and preliminary results. We refer  to the textbooks  \cite{Y}    and   to the excellent   paper  \cite{Ca} for a   survey on bases and frames in Hilbert spaces.

\subsection{Frames and  Riesz bases}
  
 A sequence of vectors ${\mathcal B}=\{v_1 ,v_2 ,v_3 ,.....\}$  in a Hilbert space  $H$ is a {\it frame} if 
there exist  constants $A, \ B>0$    such that for every $w\in H$,
\begin{equation}\label{e-frame}
 A||w||^2\leq  \sum_{j=1}^\infty |<w, v_j>|^2\leq B ||w||^2.
\end{equation}
Here, $||\ ||=||\ ||_H$ and $\langle\ \rangle=\langle\ \rangle_H$ are the norm  and the inner product in $H$.    ${\mathcal B}$ is a   {\it  tight frame } if $A=B$.

 We say that  ${\mathcal B} $    is  a {\it Riesz basis}  if   it   is 
  %
  the image of an orthonormal  basis   (called the {\it dual basis of} {\cal B}) through a  linear, bounded and invertible operator on $H$.
Equivalently,   ${\mathcal B} $   is  a {\it Riesz basis} if it is an exact  frame, i.e, if    it ceases to be a frame when any   of
its elements is removed.  See e.g. \cite{Y} for other equivalent characterizations of Riesz basis.
 
 
 Frames are over-complete Riesz basis and provide robust, basis-like   representations of vectors in a
Hilbert space.  
Because  frames  are not necessarily linearly independent,  
they are often  more easily constructible than bases. 
For example,  assume that $D\subset P=[-\pi, \pi]\times [-\pi,\pi]$. An exponential basis of $L^2(P)$ is ${\mathcal B}=\{ e^{   i (n x+ m y)}\}_{n,\,m\in\Z}$.
Let $f\in L^2(D)$, and let  $\t f(x)=\begin{cases}  f(x) &\ if\ x\in D \cr 0 & if \ x\in P-D \cr\end{cases}$.  By the Plancherel identity, 
$  
 \dsize  \sum_{n, m\in\Z} \langle \t f, \ e^{   i (n x+ m y)}  \rangle ^2_{L^2(P)} $ $=4\pi^2   ||\t f||^2_{L^2(P)} .
$  
%
Clearly, also 
  $  \dsize 
 \sum_{n, m\in\Z} \langle  f, \ e^{ i (n x+ m y)}\rangle ^2_{L^2(D)}$ $ = 4\pi^2   ||  f||^2_{L^2(D)} $, 
and so  $ {{\mathcal B}_\vert}_D $,   the  restrictions of the functions of ${\mathcal B}$ to $D$, is a  tight frame of $L^2(D)$.  

\medskip
The argument that we have used in this example  proves  the following

 \begin{Prop}\label{P-supset}
Let $P$ be a domain in $\R^d$, and let ${\mathcal F}=\{v_n\}_{n\in\Z}$ be a frame on $L^2(P)$.  
 For every domain $D\subset P$, the set $ {{\mathcal F}_\vert}_D $  is  a frame on $L^2(D)$.  
\end{Prop}

The previous example shows  that we can always construct exponential frames 
on  $L^2(D)$, but  it is often    difficult, and sometimes impossible,  to  extract    exponential   bases   from them.  That may depend on the shape of $D$, but also on the frame itself.   K. Seip proved in \cite{S1}
 that  there are exponential  frames   of $L^2(-\pi,\pi)$    which do not contain  an exponential   basis.   
 However,   it is proved in \cite{S2} that for every $L>1$, the set  $ {\mathcal B}=\{ e^{i  \pi x  \frac n L }\}_{n\in\Z}$  (which is a frame of $L^2(-1, 1)$ by Proposition \ref{P-supset}) 
 contains an exponential basis of $L^2(-1, \,1)$; that is,  there exists a sequence of integers $\{n_k\}_{k\in\Z}$   such that 
${\cal B}'=\{e^{i\pi x \frac{n_k}L   }\}_{k\in\Z}$  is an exponential basis of $L^2(-1, 1)$.    

 The main theorem in  \cite{KN} shows that  
 ${\cal B}$    contains  a Riesz basis of $L^2(I)$  also  when  $I\subset (-L,L)$ is a multi-intervals:

\begin{Thm}\label{T-Lev} (Kozma and Nitzan) Let $I\subset\R$ be a finite union of intervals. Then there exists a sequence 
 $\{\lambda_k\}_{k\in\Z} \subset\R$ such that the functions $\{e^{2\pi i\lambda_k x}\}_{k\in\Z}$ 
 form a Riesz basis in $L^2(I)$. Moreover, if $I\subset [-L, L]$, then we can chose $\lambda_k= \frac { n_k}{2L}$, with $n_k\in\Z$.
 \end{Thm}

\subsection{Stability of Riesz bases}

Riesz bases are {\it stable}, in the sense that a small perturbation of a Riesz basis produces a Riesz basis. Let us recall   Paley-Wiener stability theorem, and the celebrated Kadec   $ \frac 14$-Theorem for  exponential bases of  $L^2(-a,a)$. The proof of both theorems can be found e.g. in \cite{Y}. Kadec  theorem was originally proved in   \cite{Ka}.

   \begin{Thm}\label{P-W} (Paley-Wiener) Let $ \{v_n\}_{n\in\N}$ be a Riesz basis for a Hilbert space $H$. Suppose that $\{w_n\}_{n\in\N}$ is a sequence of elements of $H$ such that 
   \BEL{E-pw}
  \left\Vert \sum_{j=1}^n c_j( v_j-w_j)\right\Vert \leq \lambda \left\Vert \sum_{j=1}^n c_jv_j\right\Vert
   \end{equation}
  for some constant $0<\lambda<1$  and all choices of scalars $c_1$, ... $c_n$. Then $\{w_n\}_{n\in\N}$ is a Riesz basis for $H$. 
   \end{Thm}

 Note that  if  $\{v_n\}_{n\in\N}$  is an orthonormal basis of $H$, and  $\sum_{j=1}^n c_j^2 = 1$, then the right hand side of \eqref{E-pw}
 equals to $\lambda$.

\begin{Thm}\label{T-Kadec} (Kadec)  Let  $a>0$, and let $\{\alpha_n\}_{n\in\Z}\subset \R$ be such that
\BEL{E-kad}  \left |\frac{a}{\alpha_n} -  1\right|  =\epsilon_n <\frac 1{4|n| } .\EE
  Then, ${\cal B}=\{ e^{i \pi \frac{ n}{\alpha_n} x}\}_{n\in\Z}$    is a Riesz basis for $L^2(-a,a)$.  The constant $\frac 1 {4 }$ cannot be replaced by any larger number.
  \end{Thm}
  
Kadec   theorem  has been extensively  celebrated and generalized (see e.g. \cite{A}, \cite{SZ1}).  
 Theorem \ref{T-Main11} is a generalization of  Theorem \ref{T-Kadec}. 
    Indeed, 
   the inequality \eqref{E-cond11} 
    reduces to  \eqref{E-kad} if $f(y)\equiv a$  (i.e, if $T$ is a rectangle) and $ g_n(y)\equiv \alpha_n$. 
  See also Theorem  1.1 in \cite{SZ}.

\section{ Bases  of exponentials of  $L^2$ of   trapezoids.}
\setcounter{Thm}{0} \setcounter{equation}{0}

 In this  section we prove  Theorem \ref{T-Main11}  and  its corollary, Theorem \ref{T-Stepf}.     
   We start with an    easy  construction of a Riesz basis for $L^2(T)$.
 
\begin{Lemma}\label{L-basis1}    Let $ T$ be  the trapezoid  bounded by $f$. Then, 
%
\begin{equation}\label{e-b1}  {\mathcal B} =\left \{  e^{\pi i(\frac{  n x }{f(y)} +   2  ky ) }\right\}_{n, k\in\Z} \end{equation}  is a  Riesz  basis on $L^2(T)$. 
\end{Lemma}

\noindent
{\it Proof.} Recall  that    bounded, linear and invertible  functions between     Hilbert spaces maps Riesz bases into Riesz bases. 
We let $R= [-1,1]\times [0,1]$, and $\L:  L^2(R)\to L^2(T)$  be such that
\begin{equation}
\label{e-isometry}
\L(\psi)(x,y) =  (f(y)) ^{-\frac 12}\,  \psi(x/ f(y),\, y). 
\end{equation}
$\L$ is an invertible  isometry. 
Since  $ {\mathcal C} =\left \{\frac{1}{\sqrt 2} e^{ \pi i (n x   +    2  ky)  }\right\}_{n,k\in\Z}$ is an orthonormal basis of $L^2(R)$, the set 
$ 
{\mathcal B}_1= \L ( {\mathcal C} ) =\left \{(2 f(y) )^{-\frac 12} e^{ \pi i (\frac{n x }{f(y)} +   2  ky)  }\right\}_{n, k\in\Z}$ 
is an orthonormal basis for $L^2(T)$. 
 
The map $G (\psi) (x,y)= (2 f(y))^{\frac 12} \psi(x,y)  $ is linear and invertible on $ L^2(T) $ because  $f(y)>\ell>0$, and so   ${\mathcal B} $,   the image of ${\mathcal B}_1$ through $G$, is a Riesz basis for $L^2(T)$. 
   $\Box$

 \bigskip
 \noindent
 {\it Proof of Theorem \ref{T-Main11}.}    
 By Lemma \ref{L-basis1},   ${\cal B}_1=\left\{ (2  f(y))^{-\frac 12}e^{\pi i(\frac{ n x }{f(y)} + 2 ky) }\right\}_{n,k \in\Z}$  is an orthonormal Riesz basis for $L^2(T)$.  Let $\{g_n\}_{n\in\Z}\subset L^\infty(0,1)$ be  as in \eqref{E-cond11}.
Let $\dsize \tilde {\mathcal B} =\left \{(2 f(y) )^{-\frac 12} e^{\pi i(  \frac{n  x }{g_n(y)} +    2  ky) }\right\}_{n,k\in\Z}$. We show that  
\begin{equation}\label{e-PW}\int_T \left|\sum_{n,k\in\Z}\frac{ c_{n,k}}{(2f(y))^{\frac 12}} \left(e^{\pi i(  \frac{n  x }{f(y)} +    2  ky)}-e^{\pi i(  \frac{n  x }{g_n(y)} +    2  ky)}  \right)\right\vert^2dxdy \leq \lambda^2   <1\end{equation}
 whenever $\sum_{n,k\in\Z}\vert c_{n,k} \vert^2 \leq 1$. By Paley-Wiener theorem (Theorem \ref{P-W}),   $\tilde {\mathcal B} $ is a Riesz basis of $L^2(T)$, and so also     $  {\mathcal B} =\left \{  e^{\pi i(  \frac{n  x }{g_n(y)} +    2  ky) }\right\}_{n,k\in\Z}$
is   a Riesz basis for $L^2(T)$. 

In order to simplify the  proof of  \eqref{e-PW},  we use the change of variables $(x, y)\to  (x f(y), \ y)$ in the  integral. With this change of variables,  $T$ is mapped into   $R=[-1, 1]\times[0,1] $, and \eqref{e-PW} reduces to 
\begin{equation}\label{e-finproof}
 \left\Vert \sum_{n,k\in\Z}\frac{ c_{n,k}}{\sqrt 2} \left(e^{\pi i(  n  x   +    2  ky)}-e^{\pi i(  \frac{n  x f(y) }{g_n(y)} +    2  ky)}  \right)\right\Vert_{L^2(R)} \leq \lambda  <1.
\end{equation}
In the rest of the proof,  $||\ ||= ||\ ||_{L^2(R)}$. We argue as in the proof of Kadec   theorem. We let 
$$ e^{\pi i ( n x  +2 ky)} - e^{\pi i ( \frac{n x f(y)}{g_n(y)} +2 ky)}    = 
 e^{\pi i ( n x  +2 ky) }  \left(1 - e^{i \delta_n  x}\right)  $$
where $\delta_n= \delta_n(y)=\pi  n \left(\frac {f(y)}{g_n(y)} - 1\right)$.  By \eqref{E-cond11}
\begin{equation}\label{e-L} L=\sup_{{y\in[0,1]}\atop {n\in\Z}}  |\delta_n| = \sup_{n\in\Z}\left(\pi |n| \sup_{y\in[0,1]}\left| \frac{f(y)}{g_n(y)}-1\right|   \right)   <  \frac \pi {4} .
\end{equation}
We 
expand the function $1- e^{i \delta_n x}$  in $L^2(-1, \ 1)$  as a Fourier series relative to the complete orthonormal system  $\left\{1, \cos \left(  x\pi n \right) , \sin  \left( (n - \frac{1}{2})\pi x \right) \right \}_{n=1}^\infty$. We  obtain
$$1 - e^{i\delta_n x} = \left(1 - \frac{\sin  \delta_n }{ \delta_n}\right) + 2\sum_{m=1}^\infty \frac{(-1)^m   \delta_n \sin   \delta_n }{m^2 \pi^2 -  \delta_n^2} \cos  \left( \pi m x  \right) $$ 
$$+ 2i \sum_{m=1}^\infty \frac{(-1)^m   \delta_n \cos   \delta_n }{(m - \frac{1}{2})^2 \pi^2 -   \delta_n^2} \sin \left(  \left(m - \frac{1}{2}\right)\pi x  \right).$$ 
To estimate  the $L^2$  norm in \eqref{e-PW}, 
we use the above Fourier series, we change the order of the summation and we use the triangle inequality; thus,
\begin{eqnarray*}
 &\ &  \left\Vert \sum_{n,k\in\Z} \frac{c_{n,k}}{\sqrt 2}  e^{\pi i(   n x +2ky)}\left(1 - e^{i \delta_n  x}\right)  \right\Vert  
 \\ &\leq& \frac{1}{\sqrt 2}\left\Vert\sum_{n,k\in\Z} \left(1-\frac{\sin \delta_n }{ \delta_n}\right) c_{n,k}e^{ \pi i(  n  x  +2   ky) }  \right\Vert
 \\ &+ &\sqrt 2\sum_{m=1}^\infty \left\Vert\cos  \left(  m\pi x  \right) \sum_{n,k\in\Z} \frac{(-1)^m  \delta_n \sin \delta_n }{m^2 \pi^2 -  \delta_n^2}  c_{n,k}e^{ \pi i(  n  x  +2   ky) }   \right\Vert  
\\ &+ &\sqrt 2\sum_{m=1}^\infty \left\Vert
\mbox{$\sin  \left( \left(m - \frac{1}{2}\right)\pi x \right)$} \sum_{n,k\in\Z} \frac{(-1)^m    \delta_n \cos  \delta_n }{(m - \frac{1}{2})^2 \pi^2 -   \delta_n^2} c_{n,k}e^{ \pi i(   n  x   +2   ky)   }  \right\Vert 
\\ &= &  A+B+C .\end{eqnarray*}
%
%
\noindent
\eqref{e-finproof} follows if we prove that $A+B+C<1$;  
to estimate $A$, 
 we use  Plancherel theorem.
\begin{eqnarray*}
 A^2   
 &= &\frac 12 \int_0^1   \int_{-1}^{1}\left|\sum_{k \in\Z} \sum_{n \in\Z} \left(1-\frac{\sin \delta_n }{ \delta_n}\right)c_{n,k} e^{in  \pi x  + 2\pi i ky}\right|^2 dx dy
\end{eqnarray*}
\begin{eqnarray*}
& =&
\frac 12\int_{0}^{1}  \left(\int_{-1}^1  \left|\sum_{n \in\Z} \left(1-\frac{\sin \delta_n }{  \delta_n}\right) \left(\sum_{k \in\Z} c_{n,k}  e^{ 2\pi i ky} \right) e^{in  \pi x  }\right|^2 dx\right) dy
\\ &
 = &
\int_{0}^{1}   \sum_{n \in\Z} \left(1-\frac{\sin \delta_n }{ \delta_n}\right)^2 \left| \sum_{k \in\Z} c_{n,k}  e^{ 2\pi i ky} \right|^2   dy
\end{eqnarray*}
  and since  the function $t\to 1-\frac{\sin t}{t}$ is increasing  when $t\in[0,\,\pi]$, and  $|\delta_n| <L< \frac \pi 4$ by \eqref{e-L},  we 
 obtain 
\begin{eqnarray*}
 A^2 &\leq&  \left(1-\frac{\sin L}{L}\right)^2  \sum_{n \in\Z} \int_{0}^{1}  \left| \sum_{k \in\Z} c_{n,k}  e^{ 2\pi i ky} \right|^2   dy
\\
&=&\left(1-\frac{\sin L}{L}\right)^2\sum_{k \in\Z} c_{n,k} ^2 = \left(1-\frac{\sin L}{L}\right)^2
 .\end{eqnarray*}
   We argue in a similar way to show that 
 $$B \leq  \sum_{m=1}^\infty \frac{  2\, L  \sin  L}{\pi^2 m^2 -    L^2  }, \quad 
 C
  \leq \sum_{m=1}^\infty \frac{   2\,L  \cos  L}{\pi^2  \left(m -\frac{1}{2}\right)^2 -   L^2  }. $$
Thus,
$$A+B+C \leq 1-\frac{\sin L}{L}+ \sum_{m=1}^\infty \frac{\frac{2 L}{\pi} \sin  L}{\pi \left(m^2 - \frac{  L^2}{\pi^2}\right)} + 
\sum_{m=1}^\infty \frac{\frac{2 L}{\pi} \cos  L}{\pi \left(\left(m -\frac{1}{2}\right)^2 - \frac{  L^2}{\pi^2}\right)}. $$
The series $\sum_{m=1}^{\infty} \frac{\frac{2 L}{\pi}}{\pi \left(m^2 - \frac{ L^2}{\pi^2}\right)}$ and $\sum_{m=1}^{\infty} \frac{\frac{2 L}{\pi}}{\pi \left(\left(m -\frac{1}{2}\right)^2 - \frac{  L^2}{\pi^2}\right)}$ are the partial fraction expansions of the function $ \frac{1}{ L} - \cot  L$ and $\tan  L$, respectively. Hence, 
 $\lambda=A+B+C=  1 - \cos  L + \sin  L .$ 
 Since  $ L  <\frac{\pi}{4}$, we have  $\lambda <1$, and \eqref{e-finproof} is proved.

To show that the constant in \eqref{E-cond11} cannot be replaced by a smaller constant, we use a straightforward  generalization of an  example by Ingham.  
We let   $g_n(y)= \frac{ n f(y)}{n-\frac{ \rm{sgn}( n)} {4 }} $ if $n\ne 0$, and $g_0(y)=0$. Here sgn$(n)= \frac n{|n|}$.
The $g_n$'s satisfy \eqref{E-cond11}  with $\delta_n=\frac 1{4n}. $
   Let $ v_{n,k}(x,y)= e^{\pi i (\frac{nx}{g_n(y)}+ 2ky)}$. If $\{v_{n,k}\} _{m\in\Z}  $  was  a Riesz basis of $L^2(T)$, then $\{v_{n,k}(x, \ov y)\}_{n\in\Z}$ would be a Riesz basis of  $L^2(- f(\ov y), \ f(\ov y)) $ for a.e.   $\ov y\in (0,1)$.   By a change of variables,  the set  $\left\{e^{  i\pi    x(n-\frac{\rm{ sgn}( n)} 4 )}\right\}_{n\in\Z}$ would be   a Riesz  basis  of  $L^2(-1, \ 1)$, but  in \cite{I}  is proved that this is not the case.
 $\Box$

   \medskip\noindent
   {\it Proof of  Theorem \ref{T-Stepf}}.  
Assume that $f$ is continuous  (the proof of the general case is very similar).  Without loss of generality,      $ \sup_{y\in [0,1]} |f(y)|=1$. 
We construct a family of step functions $\{s_n\}_{n\ge 1}$ that satisfy the assumptions of    Theorem \ref{T-Main11};   
since we have assumed $|f(y)|\leq 1$,  \eqref{E-cond11} follows  if we prove that
\BEL{e-cond1} \sup_{y\in[0,1]}\left| \frac{1}{s_n(y)}-\frac 1{f(y)}\right| < \frac 1{4n}.
\EE
Since $f>\ell>0$,   $\frac{1}{f(y)}$ is continuous, and thus also    uniformly continuous  in $[0,1]$. If we let $\epsilon=\frac{1}{4n}$, there exists $\eta=\eta(n)>0$ such that $\left| \frac{1}{f(z)}-\frac 1{f(y)}\right|<\frac{1}{4n}$ whenever  $|z-y|<\eta $. We   partition $[0,1]$ with $N= N(n)$ intervals of    length $\frac 1N<  \eta $. For $1\leq j\leq N$, we    let $y_j= \frac {j-1}N$;  we can see at once that 
$$\left| \frac{1}{f(\frac {j-1}N)}-\frac 1{f(y)}\right|<\frac{1}{4n} \quad \mbox{ if\quad }  \frac {j-1} N< y< \frac{j}N . 
$$ 
We  let $s_n(y) \equiv   f(\frac{j-1}{N})$ if $\frac {j-1} N\leq y< \frac{j}N$ for every $1\leq j\leq N-1$, and 
$s_n(y) \equiv   f(\frac{N-1}{N})$ if $\frac {N-1} N\leq y\leq 1$.  The  $s_n$'s are step functions that satisfy \eqref{e-cond1}.   Since 
$$ \sup_{y\in [0,1]}\left|  {s_n(y)}- {f(y)}\right| =   \sup_{y\in [0,1]}s_n(y) f(y)\left| \frac{1}{s_n(y)}-\frac 1{f(y)}\right| < \frac 1{4n}, 
$$ 
 the $s_n$'s converge uniformly to $f$, as required.  $\Box$

\section {Proof of  Theorem \ref{T-Lev2}} 
\setcounter{Thm}{0} \setcounter{equation}{0}

  Let $R\subset [-1,1]\times [0,1]$ be a multi-rectangle with $N$ steps in $\R^2$.  We  let $h=  \frac 1N$, and $R= \cup_{j=1}^{N} \overline{ R_j}$, where 
  $R_j=
  (-b_j, b_j)\times ((j-1)h,\ jh )$. Recall that we have assumed  
 $0<b_j \leq 1$;  
after a dilation, we  can assume   $h=1$.  

To prove  Theorem \ref{T-Lev2} we  associate to $R$ a   multi-interval $  I$, and we define an isometry $\L: L^2 (\overline I\times [0,1])\to L^2 (R)$; 
  then we construct  an exponential basis 
  of $L^2(\overline I\times [0,1])$, and we show that 
it is mapped by $\L$ into an  exponential basis of $L^2(R)$.

Let $\vec v_j= ( 2(j-1),\ -(j-1)  )$,   so that      
$$ \tau_{\vec v_j}   R_j= (-b_j +2(j-1),\
 b_j + 2(j-1)) \times (0, \ 1) $$  
and  $\dsize \cup_{j=1}^N \tau_{\vec v_j} R_j=  I\times (0,1)$, where we have let
\begin{equation}\label{e-multis}I= \cup_{j=1}^{N } (-b_j +2(j-1),\ b_j +2(j-1)).  
\end{equation}
 Note that $I\subset (-b_1,\ 2N-1)\subset (-1,\ 2N-1)$.
 
 The segments in \eqref{e-multis}  do not overlap: indeed, for every $j\ge 1$, \newline
 $  b_j  +2(j-1) \leq -b_{j+1} +2j \iff    b_j + b_{j+1}  \leq 2  
 $   which is true because  by assumption $  b_j  \leq 1$.    

 \medskip
Let $\L: L^2 (\overline I\times [0,1])\to L^2 (R)$
 $$\L \psi (x,y)= \sum_{j=1}^{N } \chi_{\overline R_j}(x,y)\,\psi\circ \tau_{\vec v_j}(x,y).  
 $$
%
It is easy to verify that $\L$ is a linear  invertible  isometry, and hence it maps Riesz basis into Riesz basis.

Let    ${\mathcal B}= \{   e^{   2\pi i \lambda_k x}  \}_{k\in\Z}$
be an exponential basis of $L^2(I)$. By  Theorem \ref{T-Lev}, we can  assume  $\lambda_k= \frac{n_k}{2N}$ where $n_k\in\Z$.
For a fixed $k\in\Z$, we chose   $\{  m_h\}_{h\in\Z}=  \{ m_{k(h)}\}_{h\in\Z} \subset \Z $ with the following properties: the sequence $\left\{e^{  2\pi i m_h y  }\right\}_{h\in\Z}$ is a Riesz basis of $L^2(0,1)$, and 
\begin{align}\nonumber
\L   \left( e^{  2\pi i(y m_h+ x\lambda_k) }\right) &= \sum_{j=1}^{N } \chi_{\overline R_j}(x,y) \dsize e^{  2\pi i ( \lambda_k (x+2(j-1) )+m_h (y-j+1) )  } 
\\ \label{e-inv} &=
 \chi_{R}(x,y) e^{  2\pi i(  x\lambda_k+y m_h) }.
\end{align}
We have \eqref{e-inv}  if   $2\pi(j-1) (2\lambda_k  -m_h ) = 2\pi(j-1) ( \frac{n_k}{N}  -m_h ) $ is an integer multiple of $2\pi $. If we   let
$m_h=  \{n_k\} +   h $,
 where  $\{  n_k\}$ is  the remainder  of the division of $ n_k$ by $N$,   
 the sequence  $\left\{e^{   2\pi i m_h y  }\right\}_{h\in\Z}=   \left\{e^{ 2\pi iy \left(h+   \{n_k\} \right)  }\right\}_{h \in\Z}$  is an exponential  basis of $L^2(0,1)$,  and  \eqref{e-inv} is satisfied.
By Lemma  \ref{L-multibase} below, the set   $\B_1= \left\{e^{   \pi i(y m_h+ x\lambda_k)  }\right\}_{h,k\in\Z}$ is an exponential basis of \newline $L^2(\overline I\times [0,1])$, and so $\L (\B_1)$ is an exponential  basis of $L^2(R)$. $\Box$

 \medskip

We are left to  prove the following

\begin{Lemma}\label{L-multibase}
Let $D$  be a domain of $\R^k$ and   $E$ a  domains  of $\R^d$. Let   $\{v_n(x)\}_{n\in\Z}$ be a  Riesz basis of  $L^2(D)$;  assume that  for every $n\in\Z$,  there exists a sequence  $\{n(m)\}_{m\in\Z}\subset\Z$ and a Riesz basis $ \{ w_{n(m)}(y)\}_{m\in\Z}$   of $L^2(E)$.
Then ${\mathcal B}=\{v_n(x) w_{n(m)}(y)\}_{n,m \in\Z}$ is Riesz basis of $L^2(D\times E)$.
\end{Lemma}

\medskip
\noindent
{\it Proof.}  
Let  $\{v'_n(x)\}_{n\in\Z}$  and  $\{w'_{n(m)}(y)\}_{m\in\Z}$) be   dual bases of $\{v_n(x)\} $    and $\{w_{n(m)}(y)\}$.  Let $f(x,\, y)\in L^2(D\times E)$.  Then, for a.e.  $y\in E$, the function  $x\to f(x,y)$ is in $ L^2(D)$, and   
$$ f(x,y)= \sum_{n=-\infty}^\infty  \langle f,\ v'_{n }\rangle_{D}  v _{n }(x).$$
Let $f_n(y)= \langle f,\ v'_n \rangle_D$. It is easy to verify that   $f_n \in L^2(E)$,  and  so  
$$   f_n(y)=  \sum_{m=-\infty}^\infty \langle f_n,\ w'_{n(m)}\rangle_{E}  w _{n(m)}(y) . $$
By Fubini  theorem, $ \langle f_n, \   w'_{n(m)}\rangle_{E} =$   $\langle   \langle f, \ v'_n\rangle_D,\ w'_{n(m)}\rangle_{E}= \langle f, \ v'_n w'_{n(m)}\rangle_{ D\times E }$; therefore,
\begin{align*}
f(x,y)& = \sum_{n=-\infty}^\infty f_n(y)v_n(x) =  \sum_{n=-\infty}^\infty\left(\sum_{m=-\infty}^\infty  
\langle f_n,\ w'_{n(m)}\rangle_{E}\, w_{n(m)}(y)\right) v_n(x)
\\ 
&=\sum_{n,m=-\infty}^\infty\langle f, \ w'_{n(m)}v'_n\rangle_{ D\times E }\  w_{n(m)}(y)v_n(x)
\end{align*}
which shows that  $w_{n(m)}(y)v_n(x)$ is a  Riesz basis of $L^2(D\times E)$  with dual  basis $w'_{n(m)}(y)v'_n(x)$.
 $\Box$

\section{Spherical trapezoids}
\setcounter{Thm}{0} \setcounter{equation}{0}

In this section we prove a multi-dimensional   version of Lemma \ref{L-basis1}  and we discuss generalizations of our main theorems in dimension $d \ge 2$. We will write $x\in\R^d$ as $(x', y)$, with $y\in \R$ and $x'\in\R^{d-1}$.
Let $f:[0,1]\to\R$  be as in \eqref{e-admis}.  The  {\it trapezoid with spherical basis bounded by $f$} is the set
\begin{equation}\label{e-Trapezd} T=\{ (x',y)\in\R^{d} : \ |x'| \leq f(y),  \  0\leq y\leq 1 \}.  \end{equation} 
 We say that $T$  is a  {\it  multi-cylinder }  if $f$ is a regular step function.

We let  $L^2_S(T)$ be the set of $ L^2(T)$  functions    which are radial in $x'$.  
When $d=2$,    a  multi-cylinder  is   a multi-rectangle and $L^2_S(T)$ is the set of all $L^2(T)$ functions   that are even in $x'$. 

The following lemma reduces the construction of a Riesz basis of $L^2_S(T)$ to that  of a Riesz basis of $L^2(R)$, where $R=[0,1]\times [0,1]$.
It is the analogue of Lemma \ref{L-basis1}.
  
  \begin{Lemma}\label{L-basisd}    Let $ T $ be the  spherical trapezoid bounded by $f$. Then,  
  $$   {\mathcal B} =\left \{  e^{2\pi i (\frac{n |x'| }{f(y)} +    ky)  }\right\}_{n, k\in\Z}$$ is a  Riesz  basis of  $L^2_S(T)$. 
\end{Lemma}

\noindent
{\it Proof.}  
Let   $\L:L^2(R)\to L^2_S(T)$  be such that
$$
\L(\psi)(x',y) = (|S^{d-2}| \ f(y))^{-\frac 12} r^{-\frac{d-2}{2}} \psi(r/f(y),\, y)   
$$
where  $r=|x'|$ and  $|S^{d-2}|$ is the measure of the unit sphere in $\R^{d-1}$. When $d=2$, we let $ |S^{0}| =1$.
  $\L$ is an invertible  isometry. Since   $\dsize {\mathcal C} =\left \{  e^{ 2\pi i (n r   +    ky)}\right\}_{n,k\in\Z}$ is an orthonormal basis of $L^2(R)$, the set 
$$
{\mathcal B}_1= \L( {\mathcal C} )=\left \{(|S^{d-2}|\ f(y))^{-\frac 12} r^{-\frac{d-2}{2}} e^{ 2\pi i (\frac{n r}{f(y) }  +    ky)} \right\}_{n, k\in\Z}$$
is an orthonormal basis of $L^2_S(T)$. We  conclude that  $\B$ is a Riesz basis of $L^2_S(T)$  as in  Lemma  \ref{L-basis1}.
   $\Box$
 
 \medskip
We can prove versions of  Theorems 
  \ref{T-Stepf}, \ref{T-Main11} and    \ref{T-Lev2}  for spherical trapezoids.  
The analogue of Theorem \ref{T-Stepf} is
  
 \begin{Thm}\label{T-Main111}  Let 
  $T\subset \R^d$ be the spherical trapezoid bounded by a piecewise continuous $f$. There exists a family of regular step functions 
$\{s_n \}_{n\ge 1} $  which converge  uniformly   to $f$, for which 
 ${\mathcal B} =\left \{   e^{ 2\pi i ( \frac{n |x'| }{s_n(y)} +    ky) }\right\}_{n,k\in\Z} $
is a Riesz basis of $L^2_S(T)$.  
 \end{Thm}

Theorem \ref{T-Stepf} is not a special case of Theorem \ref{T-Main111}: indeed, when   $d=2$, the set $\B$  defined above is not  the same   $\B$   in Theorem \ref{T-Stepf}.
The proof of Theorem \ref{T-Main111} is    
almost a line-by-line repetition of the proofs of Theorems \ref{T-Stepf}; it follows from a straightforward  generalization of  Theorem \ref{T-Main11}.
 
The analogue of Theorem \ref{T-Lev2} is 
\begin{Thm}\label{T-Lev3}
Let    $T\subset B_1(0)\times [0,1]$ be  a multi-cylinder with $N$ steps. There exists sequences  $\{n_k\}_{k\in\Z}$,  $\{m_h\}_{h\in\Z}\subset \Z$   such that  $\{e^{2\pi i (\frac{  n_k |x'|}{ N}+ N m_hy)  }\}_{h, k\in\Z}$ 
is a  Riesz basis of $L^2_S(T)$.
 \end{Thm}
 
  %
Our method of proof  allows to  construct Riesz bases, but not exponential bases of $L^2_S(T)$. The existence of exponential bases on domains with radial symmetry is still an open problem.  

 \section{Remarks and open problems} We have constructed exponential bases of $L^2(R)$ for certain   multi-rectangles $R$; we are wondering if is it possible to do the same  when  $R$ is  a union of a generic   family of disjoint rectangles.     In  Theorem \ref{T-Lev2}    the rectangles in $R$ have the same height, and this fact allows to reduce the construction of     an exponential basis of $L^2(R)$ to that of an exponential basis of $L^2(I)$.
Our construction   does not seem to work well   for  general multi-rectangles. 
 
 Complex analysis methods  have   often been used in these problems.  We cite, for example, the recent paper by J. Marzo \cite{M};    the author  proves the existence of a Riesz basis of
exponentials on a finite union of congruent cubes of $\R^n$    by finding  complete interpolating sequences in a suitable Paley-Wiener space.  We believe that    the    proof in \cite {M}  cannot be easily generalized when the cubes are replaced by parallelepiped.

 We are also wondering if, for certain trapezoids $T$,  the  construction of   exponential  Riesz bases on a family of multi-rectangles $R_n$  that approximate  $T$ (in the sense that the measure of the symmetric difference of $R_n$ and $T$ goes to zero when $n\to \infty$) can lead to the construction of an exponential basis of $L^2(T)$. We   plan to pursue this  investigation in a subsequent paper.

   \medskip
LAURA DE CARLI,
  Florida International Univ., Miami (FL) 33199
  %
  
  E-mail: decarlil@fiu.edu
  
  ANUDEEP KUMAR, 
  The George Washington Univ., 
 Washington, DC 20052. 
E-mail: anudeep@email.gwu.edu

\end{document}